\def\TNRNR{\overline{T}^*}
\def\Fbar{\overline{F}}
\def\signFzero{\Sigma_0}
\def\signFbarzero{\overline{\Sigma}}

\def\te{\theta} 
\def\tedd{{\ddot \theta}} 

\def\x{x} 
\def\xd{{\dot x}} 
\def\xdd{{\ddot x}} 

\def\y{y} 
\def\yd{{\dot y}} 
\def\ydd{{\ddot y}} 

\def\z{z} 
\def\zd{{\dot z}} 

\def\xso{{\bf x}^*(t)} 

\def\xsdo{\dot {\bf x}^*(t)} 

\def\xs{{\bf x}(t)} 
 
\def\xsd{\dot {\bf x}(t)} 

\def\ut{{\bf{u}}(t)}

\def\systemA{{\bf{A}}}
\def\systemB{{\bf{B}}}

\def\po{{\bf p}^*(t)} 
\def\pdo{\dot {\bf p}^*(t)} 
\def\p{{\bf p}(t)} 
 
\def\ue{u^e(t)} 

\def\u*{{\bf{u}^*}(t)} 

\def\utwo*{{\bf{u}_2^*}(t)} 
\def\D{\Delta} 
\def\de{\delta} 

\def\sg{\left\lbrace}

\def\sq{\left\lbrack}
\def\dq{\right\rbrack}
\def\st{\left (}
\def\dt{\right )}
\def\oftime{\!\left( t \right)}

\def\ds{\displaystyle}

\def\be{\begin{equation}}
\def\ba{\begin{array}}
\def\ee{\end{equation}}
\def\ea{\end{array}}
\def\bdm{\begin{displaymath}}
\def\edm{\end{displaymath}}

\def\bsy{\boldsymbol}

\def\Tofxi{T^*\!\left(\boldsymbol{\xi}(\tau)\right)}

\documentclass[twocolumn]{autart}    

\usepackage{booktabs} 
\usepackage{makecell} 

\usepackage{graphicx}          
\usepackage{amsmath} 
\usepackage{amssymb}
\usepackage[dvips]{epsfig}    
\usepackage{epstopdf}

\usepackage{scalerel} 

\usepackage{enumitem}
\usepackage{enumitem}
\usepackage{amsfonts}
\usepackage{multirow}
\usepackage[utf8x]{inputenc} 
\usepackage{color}

\newcommand{\rosso}[1]{\textcolor{red}{#1}}

\begin{document}

\begin{frontmatter}

\title{Analytic Solution of the Time-Optimal Control of a Double Integrator from an Arbitrary State to the State-space Origin} 
\author[NPS]{ Marcello Romano}\ead{mromano@nps.edu} , and
\author[Roma]{Fabio Curti}\ead{fabio.curti@uniroma1.it}            

\address[NPS]{Dept. of Mechanical \& Aerospace Engineering, Naval Postgraduate School, 700 Dyer Rd., Monterey, California 93940. U.S.A.}             
\address[Roma]{School of Aerospace Engineering, Sapienza University of Rome, Via Salaria, 851, 00138 Rome. Italy.}  


          
\begin{keyword}                           
Analytical design; Rest to Rest Control; Regulation              
\end{keyword}                             

\begin{abstract}                          
This brief note presents known results about the minimum-time control of a double integrator system from an arbitrary initial state to the state-space origin (minimum-time regulation problem, or special problem). The main purpose of this note is didactical. Results are presented in all details and following a step by step procedure. 
\end{abstract}

\end{frontmatter}

\section{Introduction}

The theory of minimum-time control of linear systems is for the most part maturely established~\cite[115--187]{Pontryagin:1962aa}~\cite[127--158]{Lee1967}\cite[395--426]{athans1966}\cite[83--173]{Bolt1971}\cite[248--249]{Kirk1998}~\cite[127--158]{Lee1967}.

Among linear systems, the double integrator is widely studied as it constitutes useful model for many dynamic phenomena encountered in engineering and science~\cite{Rao2001}. 

This note collects and presents in all details results that are originally found in many textbooks, including for instance~\cite[115--187]{Pontryagin:1962aa},~\cite{athans1966}. 

Results are here presented and demonstrated in a step by step fashion, that is deemed particularly useful for students. 

\section{System Dynamics}
A one-d.o.f. \emph{double-integrator} system is governed by the following equation
\begin{equation} \label{doubleintegrator}
\sg 
\ba{l}
I \ydd \oftime= C\oftime\\
|C\oftime| \le C_{\text{max}}  
\ea
\right. ,
\end{equation}
where $I\in \mathcal{R}$ is the inertia parameter, $y\oftime\in\mathcal{R}$ is the displacement, $C\oftime\in\mathcal{R}$ is the control and $C_{\text{max}} $ is the maximum magnitude of the control. Without loosing generality, the equations above can be rewritten in a convenient scaled form. In particular, by transforming the displacement variable to the new variable $x\oftime\in\mathcal{R}$ and the control variable to the new variable $u \oftime\in\mathcal{R}$ as follows
\begin{equation}
x\oftime = \displaystyle\frac{I}{C_{\text{max}}} \, y \oftime, \qquad
u\oftime  =\displaystyle \frac{C\oftime}{C_{\text{max}}},
\end{equation}
the system of  Eq.~\ref{doubleintegrator} can be equivalently written as
\begin{equation} \label{doubleintegrator2}
\sg 
\ba{l}
\xdd \oftime= u\oftime\\
|u\oftime| \le 1  
\ea
\right. ,
\end{equation}
or, in state-space form, as
\be \label{DI_statespace}
\sg 
\ba{l}
 \xsd= {\bf A}  \xs+{\bf B} u(t)\\
|u\oftime| \le 1  
\ea
\right. ,
\ee
where
\be \label{DI_AB_matrix}
 {{\bf x} \oftime}=\sq \ba{c} \x \\ \xd  \ea \dq; \hskip.1cm {\bf A} = \sq \ba{cc} 0&1 \\ 0&0 \ea \dq; \hskip.1cm {\bf B} = \sq \ba{c} 0 \\ 1 \ea \dq.
\ee
\section{Solution of the Special Problem}
\begin{figure}[h!]
		\includegraphics[width=\linewidth]{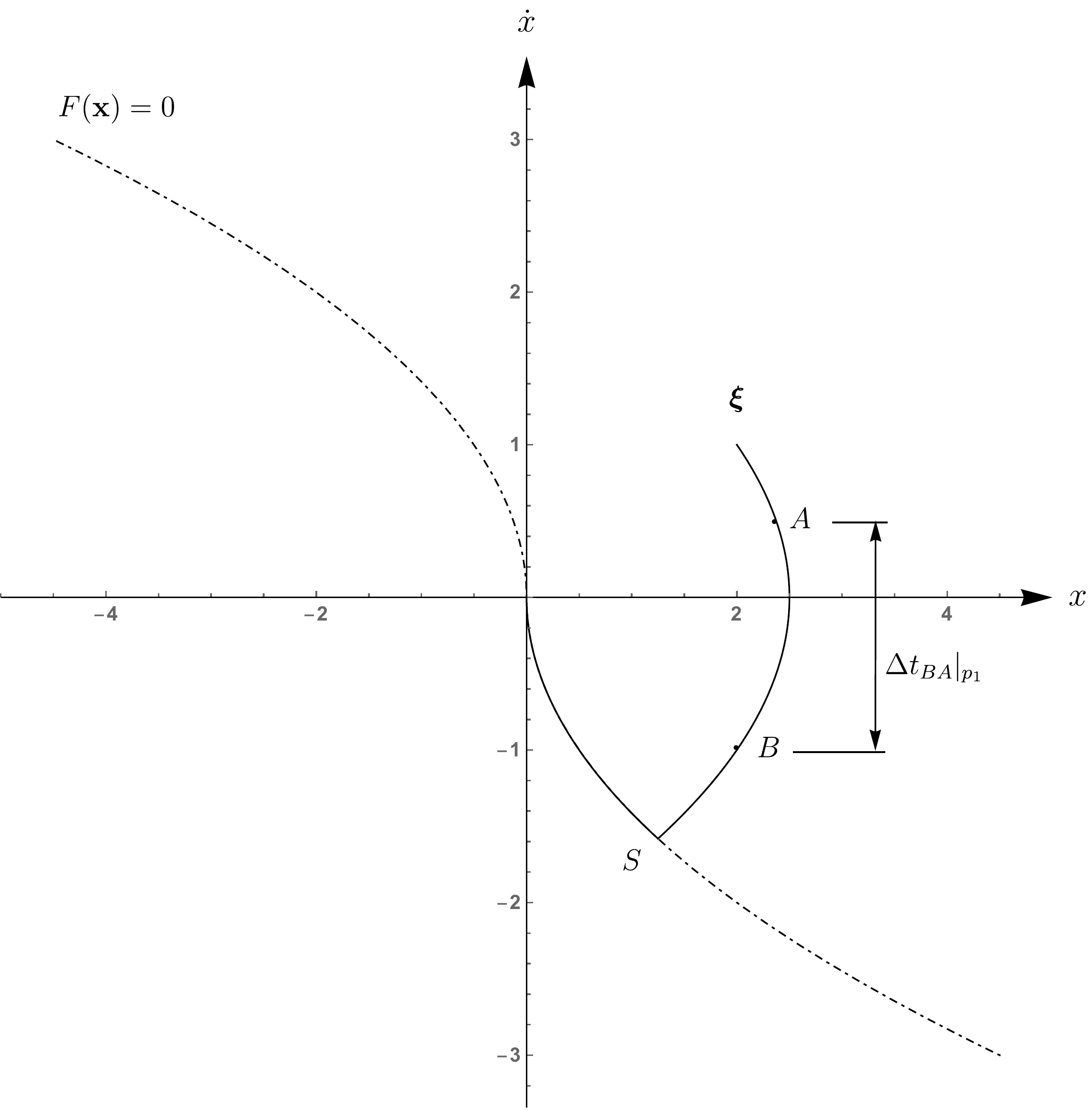}
	\caption{Double integrator: sample no-rest to rest (i.e. special problem) optimal state trajectory for an initial state such that $F_0>0$ (initial state: $\bsy \xi=[2,1]'$)}
	\label{DI_traj_nr2r}
\end{figure}

This section reports known results regarding the time-optimal control of the double integrator between an arbitrary initial state and the state-space origin. This problem can be referred as \emph{No-rest to rest control to the origin} or as \emph{Special Problem}. 
\begin{thm}\label{teoremaDINR2R}
{\textit {\textbf {(Double Integrator: Solution of the Minimum-time optimal control problem from an arbitrary state to the origin)}}}\\ 
Assume boundary states 
\be \label{DI_boundaryconditionNR2R}
 {{\bf x} (0)}=\boldsymbol{\xi}\triangleq\sq \ba{c}x_0 \\ \dot x_0  \ea \dq; \hskip.5cm {\bf{x}}(T^*)= \bf{0}=\sq \ba{c}0 \\ 0   \ea \dq.
\ee
Define
\be \label{definitionofF}
F({\bf x} )=F(\x,\xd)\triangleq   \x+ {\rm sgn} \!\left( \xd \right)\frac{ \xd^2}{2}.
\ee
In particular, it is named \emph{switching curve} the following curve composed of two arcs of semi-parabolas joined at the origin
\be \label{switchcurve_DI}
F({\bf x} )=0.
\ee
Define, furthermore,
\begin{align}
\label{DINR2RF0}
&F_0\!\!\! &&\triangleq F(\boldsymbol{\xi}) = F(\x_0,\xd_0) =x_0+ {\rm sgn} \!\left( \xd_0 \right) \frac{ \xd_0^2}{2},\\
&\signFzero\!\!\!  &&\triangleq {\rm sgn} (F_0). \label{DINR2RSigma0def}
\end{align}
The optimal solution of the problem is as follows

(a) if $F_0=0\leftrightarrow\signFzero=0$, the optimal control history is
\be \label{controlonswitchcurve}
u^* \oftime = u^* = -{\rm sgn} (\xd_0 ), \hskip.4cm \, t \in [0, T^*],
\ee
with
\be\label{mintimeDINR2R1}
T^*=-u^*\,\xd_0,
\ee
or equivalently
\be\label{DINR2RdiseqTstar0}
T^*=|\xd_0|.
\ee
Furthermore, the optimal trajectory is
\be \label{stateequationevo1}
 {{\bf x}^* \oftime}=\sq \ba{c} \x^*(t) \\ \dot{x}^*(t)  \ea \dq= \sq \ba{l} \displaystyle u^*\frac{ t^2}{2} + \xd_0 t +x_0\\u^* t +   \xd_0  \ea \dq, \quad t\in[0,T^*],
 \ee
 which corresponds in the phase-plane to the parabola
 \be \label{parabola1}
\st \x-{1 \over 2} u^* \xd^2 \dt =\st \x_0-{1 \over 2} u^* \xd^2_0 \dt.
\ee
\vskip.05cm
(b) if $F_0\neq0\leftrightarrow \signFzero=\pm1$, the optimal control history is
\begin{equation} \label{controloutofswitchcurve}
u^* \oftime =
\left\{
\begin{alignedat}{2}
u^*_1 &= -\signFzero, \quad  &t &\in [0, \Delta_1)\\
u^*_2 &= -u^*_1 = \signFzero,\quad &t &\in [\Delta_1,T^*]
\end{alignedat}
\right.,
\end{equation}
where the duration of each ``bang'' is~\footnote{In the particular case when $(x_0\neq 0,\,\xd_0=0)$ it yields $\Delta_1=\Delta_2=\Lambda_0=\sqrt{|x_0|}$.}
\be \label{mintimeDINR2R2}
\sg \ba{l} \ds
 \D_1= \Lambda_0+\signFzero\, \xd_0  \\ \ds
\D_2=\Lambda_0 
\ea \right.,
\ee
with
\be \label{DINR2Rlambdef}
\Lambda_0\triangleq\sqrt{\signFzero \,\x_0+\frac{\xd^2_0}{2} },
\ee
where
\be\label{DINR2Rlambineq} 
\signFzero \,\x_0+{1 \over 2} \xd^2_0 > 0, \quad \forall \,{\bf{x}}(0)\neq{\bf{0}},
\ee
and therefore
\be \label{DINR2Rlambdeflt0}
\Lambda_0 \in \mathcal{R}^+,  \quad \forall \,{\bf{x}}(0)\neq{\bf{0}},
\ee
Moreover, the following inequality yields (which is more restrictive than the inequality implied in Eq.~\ref{DINR2Rlambdeflt0})
\be\label{DINR2Rlambineq33} 
\Lambda_0> |\xd_0|, \quad \quad \forall \,\left(\signFzero,\xd_0 \right) :\!\!| \,\, \signFzero\, \xd_0 = -|\xd_0| <0,
\ee
which immediately yields
\be
\D_1>0.
\ee
It yields, finally, 
\be \label{DINR2RTstarexplicit}
T^*= \D_1+ \D_2= 2\Lambda_0+\signFzero\, \xd_0,
\ee
which enjoys the property
\be \label{DINR2RdiseqTstar}
T^*> |\xd_0|.
\ee
Furthermore, the optimal trajectory is
\be \label{stateequationevo2}
 {{\bf x}^* \oftime}=\sq \ba{c} \x^*(t) \\ \dot{x}^*(t)  \ea \dq= \sq \ba{l} \displaystyle u^*\frac{ t^2}{2} + \xd_0 t +x_0\\u^* t +   \xd_0  \ea \dq,  t\in[0,\Delta_1),
 \ee
 \be \label{stateequationevo22}
 {{\bf x}^* \oftime}=\sq \ba{c} \x^*(t) \\ \dot{x}^*(t)  \ea \dq= \sq \ba{l} \displaystyle u^*\frac{ t^2}{2} + \xd_s t +x_s\\u^* t +   \xd_s  \ea \dq, t\in(\Delta_1,\Delta_2],
 \ee
 where
 \be \label{intersectionpoint_DI1}
\sg \ba{l} \ds
 \x_s={1 \over 2} \left( \x_0 +{1 \over 2} \signFzero \,\xd^2_0 \right) \\ \ds
\xd_s= -\signFzero \Lambda_0
\ea \right.,
\ee
are the coordinate of the switch point $S$ (see also Fig.~\ref{DI_traj_nr2r}). 
 The optimal trajectory corresponds in the phase-plane to the union of the following two arcs of parabolas, connected at the switch point $S$.
 The first arc of parabola $p_1$ runs between the initial state and the switch point, and the parabola has equation
 \be \label{parabola1}
p_1: \quad \st \x-{1 \over 2} u^* \xd^2 \dt =\st \x_0-{1 \over 2} u^* \xd^2_0 \dt.
\ee
The second arc of parabola $p_2$ runs between the switch point and the final state (origin of the phase-plane) and the parabola has equation 
\be \label{parabola12}
p_2: \quad \st \x-{1 \over 2} u^* \xd^2 \dt =\st \x_s-{1 \over 2} u^* \xd^2_s \dt.
\ee
Finally, notably, the time elapsed between any two successive points $A$ and $B$ both on $p_1$ equates the difference in ordinate ($\dot{x}$) between the two points, i.e. (See also Fig.1)
\be \label{eq30}
\Delta t_{BA}|_{p_1}= \frac{\dot{x}_B-\dot{x}_A}{-\Sigma_0},
\ee
and analogously for any two successive points  $C$ and $D$ both on $p_2$,
\be
\Delta t_{DC}|_{p_2}= \frac{\dot{x}_D-\dot{x}_C}{\Sigma_0}.
\ee
Finally the time elapsed between any two successive points $B$ on $p_1$ and $C$ on $p_2$ is the sum of the time elapsed between $B$ and the switch point $S$ along $p_1$ and the time elapsed between $S$ and $C$ along $p_2$, i.e.
\be\label{eq32}
\Delta t_{CB}= \Delta t_{SB}|_{p_1}+\Delta t_{CS}|_{p_2}.
\ee

\end{thm}

Notably, since Theorem 1 is valid for an arbitrary initial condition, it gives the feedback optimal control synthesis for regulating to zero a double-integrator system.

{\bf Proof of Theorem \ref{teoremaDINR2R}}:\\
The system of Eq. \ref{DI_statespace} is linear time invariant, normal and has real eigenvalues. General theorems valid for this class of systems, guarantee that the time optimal control sequence exists, is unique, and has at most one control switch~\cite[th.6-5, p.399-420]{athans1966}. 

The detailed demonstration of Theorem 1 is found below by exploiting Pontryagin's principle together with geometric analysis. 

For a control history and controlled trajectory to be optimal, Pontryagin's principle requires that exists a constant multiplier $\rho^*\ge0$ and a costate $\po \in \mathcal{R}^n, t\in[0,T^*]$ such that the following conditions are safisfied~\cite[p.94 and p.108]{Heinz1}\cite[th.6-4, p.396]{athans1966}\cite[p.18]{Pontryagin:1962aa}:
\begin{enumerate}
\item{\emph{Non-triviality of the multipliers: }  
\be
\left(\rho^*,\po \right)\neq \left(0,{\bf{0}}\right).
\ee}
\item{\emph{Canonical Equations: }$\xso$ and $\po$ satisfy
\be
\ba{l}\label{problem1canonicalequations}
\xsdo= \displaystyle \frac{\partial H}{\partial {\bf{p}}} =\systemA \xso +\systemB \u*\\\\
\pdo=\displaystyle- \frac{\partial H}{\partial {\bf{x}}}=- \systemA' \po,
\ea
\ee
with boundary-state conditions as in Eq.~\ref{DI_boundaryconditionNR2R}, and Hamiltonian function given by
\be \label{eqregulatorhamiltonian}
H(\rho,{\bf{p}},{\bf{x}}, {\bf{u}}) =\rho+\p' (\systemA \xs+\systemB \ut).
\ee
The second of Eqs.~\ref{problem1canonicalequations} yields
\be \label{costateevo}
\po=\e^{-\systemA' t}{\bf\,{p^*}}(0).
\ee 
}
\item{\emph{Minimum Condition: }The Hamiltonian has an absolute minimum at $\ut=\u*$ 
\be \label{eq8}
H(\rho^*,{\bf{p}}^*,{\bf{x}}^*, {\bf{u}}^*) \le H(\rho^*,{\bf{p}}^*,{\bf{x}}^*, {\bf{u}})
\quad \forall \, t \in [0, T^*],
\ee
which requires
\be \label{regulatoroptimalcontrol}
\u*=-{\rm sgn}\!\st \systemB' {\bf{p}^*}(t) \dt, \quad \forall t\in[0,T^*].
\ee}
\item{\emph{Transversality Condition: }In case of minimum-time control and fixed boundary states the Hamiltonian is zero at the end-poin~\cite[187]{Kirk1998}. 
Furthermore, the value of the costate at the initial and final time is free.}
\item{\emph{Stationarity of the Hamiltonian: }
If conditions 2 and 3 above are satisfied~\cite[p.36]{Geering2007}, 
 by taking into account Eq.\ref{eqregulatorhamiltonian}, it yields
\be
\frac{d H}{dt}=\frac{\partial H}{\partial t}=0 \rightarrow H=constant.
\ee
Therefore, by considering condition 4, it yields
\be \label{zerohamiltonian}
H(\rho^*,\po,\xso, \u*) = 0, \,\,\, \forall \, t \in [0, T^*].
\ee
}
\end{enumerate} 
By taking into account Eqs.~\ref{DI_statespace} and~\ref{DI_AB_matrix}, and Eqs.~\ref{eqregulatorhamiltonian},~\ref{costateevo} and~\ref{regulatoroptimalcontrol}, it yields
\be 
\po=\sq \ba{c} p_x^*(t) \\ p^*_{\dot{x}}(t) \ea \dq=\sq \ba{c} p^*_{x0} \\ p^*_{{\dot{x}}0} -p^*_{x0}\,t\ea \dq,
\ee
\be \label{hamiltonianforallt}
H=\rho^*+p^*_{x0} \dot{x} +  \left(p^*_{{\dot{x}}0} -p^*_{x0}  t \right) u^*(t),
\ee
and
\be \label{regulatoroptimalcontrolexpl}
u^*(t) =-{\rm sgn}\!\st p^*_{{\dot{x}}0} -p^*_{x0}\,t \dt.
\ee
Therefore the optimal control is either $+1$ or $-1$ with at most one switch because the argument of the sign function in Eq.~\ref{regulatoroptimalcontrolexpl} is linear in the time variable, and therefore crosses zero at most once.
Moreover, the condition of optimality in Eq.~\ref{zerohamiltonian} becomes
\be \label{zerohamiltonianspecific}
H=0,\,\quad t \in [0, T^*].
\ee

By considering an interval of time $t\in[0,t]$ during which $u^*$ is constant, by integrating in time Eq.~\ref{DI_statespace}, it yields
\be \label{stateequationevo}
 {{\bf x}^* \oftime}=\sq \ba{c} \x(t) \\ \xd(t)  \ea \dq= \sq \ba{l} \displaystyle u^*\frac{ t^2}{2} + \xd_0 t +x_0\\u^* t +   \xd_0  \ea \dq,
 \ee
which, by solving the first equation for $t$ and substituting the result into the second equation and by taking into account that $u^*=1/u^*=\pm 1$, yields the following optimal control path on the phase-plane
\be \label{parabolas}
\st \x-{1 \over 2} u^* \xd^2 \dt =\st \x_0-{1 \over 2} u^* \xd^2_0 \dt.
\ee
Eq.~\ref{parabolas} represents a parabola on the phase-plane having $x$ as abscissa-axis coordinate and $\dot{x}$ as ordinate-axis coordinate. The parabola has vertex of coordinates $(\x_0-{1 \over 2} u^*  \xd^2_0, 0)$ and axis of symmetry coincident with the abscissa axis $(\xd=0)$. 
When $u^*=+1$, the parabola has positive concavity toward the positive end of the abscissa. When $u^*=-1$, the parabola has negative concavity toward the positive end of the abscissa. In both cases the parabola is run in a clock-wise fashion by the representative point of the system state on the phase plane. 

Furthermore, a curve named {\it switching curve} on the phase-plane can be defined with Eq.~\ref{switchcurve_DI}. The switching curve is the union of the two semi-parabolas passing through the phase-plane origin and occupying second and fourth quadrant (see Fig.~\ref{DI_traj_nr2r}). Each semi-parabola is an element of one of the two families of parabolas in Eq.~\ref{parabolas}.

Assume first that $F_0=0$, where $F_0$ is defined in Eq.~\ref{DINR2RF0},  i.e. assume that the initial state belongs to the switching curve. In this case, the following controlled state path equation on the phase-plane
\be \label{stateevonoswitch}
\x-u \frac{ \xd^2}{2}=0,
\ee
with the constant control history $u=-{\rm sgn} \!\left( \xd_0 \right)$, i.e., in other words, the path equation $F({\bf x} )=0$, satisfies the state equation. 
In fact, it satisfies Eqs.~\ref{stateequationevo} or the equivalent Eq.~\ref{parabolas}. Furthermore, Eq.~\ref{stateevonoswitch} satisfies the boundary conditions in Eq.~\ref{DI_boundaryconditionNR2R}. It satisfies the initial condition because, by hypothesis, $F_0=0$, and the final condition since Eq.~\ref{stateevonoswitch} becomes an identity when ${{\bf x}}={\bf{0}}$. Finally, by considering the following initial values for the costate
\be~\label{initialcostatecasef00}
p^*_{{\dot{x}}0} ={\rm sgn} \!\left( \xd_0 \right)\,\rho^*, \qquad  p^*_{x0}=0, \qquad \forall \rho^*,
\ee
the zero-value condition for the Hamiltonian in Eq.~\ref{zerohamiltonianspecific} holds true, as it can be immediately verified by considering the two possible cases of ${\rm sgn} \!\left( \xd_0 \right)=\pm1$. 
By substituting the values of initial costate variables in Eq.~\ref{initialcostatecasef00} into the argument of the sign in Eq.~\ref{regulatoroptimalcontrolexpl}, it remains confirmed that there is no switch of the control during this optimal maneuver. \newline
By observing the phase-plane geometry, it is immediate to discern that the optimal state trajectory is an arc of one of the two parabolas belonging to the switching curve. In particular it is the arc between the point $(\x_0,\xd_0)$ and the origin. Therefore, the optimal control history is as reported in Eq.~\ref{controlonswitchcurve}, with $T^*$ given by Eqs.~\ref{mintimeDINR2R1}, which is obtained by substituting the value of $\dot x(T^*)$ from Eq.~\ref{DI_boundaryconditionNR2R} into the second of Eqs.~\ref{stateequationevo}. 

Assume now that $F_0\neq 0$. By observing the phase-plane geometry and considering that the control can only be ``bang-bang'' and there can only be one switch, it is immediate to discern that the optimal state trajectory is the union of two arcs of parabolas. The first arc of parabola goes from the point $(\x_0,\xd_0)$ to the point $S$ where the switching curve is intersected, the second one, which belongs to the switching curve, goes from the point $S$ to the origin (see Fig. \ref{DI_traj_nr2r} ). If $F>0$, i.e. if the initial state is above the switching curve, the first control ``bang'' is $-1$, the second one is $+1$ and the switch point is on the arc of the switching curve in the fourth quadrant of the phase-plane. Viceversa, if $F>0$, i.e. if the initial state is below the switching curve, the first control ``bang'' is $+1$, the second one is $-1$ and the switch point is on the arc of the switching curve in the second quadrant of the phase-plane.Therefore, the optimal control history is as reported in Eq.~\ref{controloutofswitchcurve}. 
The switch point $S$ has coordinates
\be \label{intersectionpoint_DI}
\sg \ba{l} \ds
 \x_s={1 \over 2} \left( \x_0 +{1 \over 2} \signFzero \,\xd^2_0 \right) \\ \ds
\xd_s= -\signFzero \Lambda_0
\ea \right.,
\ee
as obtained by solving the algebraic system of Eq.~\ref{parabolas}  and Eq.~\ref{switchcurve_DI}, after substituting $x$ with $x_s$, $\dot{x}$ with $\dot{x}_s$, $u^*$ with $u^*_1$ from Eq.~\ref{controloutofswitchcurve}, and ${\rm sgn} \!\left( \dot{x}_s \right)$ with $-\signFzero$. This latter expression is also used to make explicit the discrimination between $+1$ and $-1$ in the equation $\xd_s=\pm \Lambda_0$, resulting from solving the algebraic system.

The duration of the control ``bangs'' is as in Eq.~\ref{mintimeDINR2R2}.  In particular, the expression of $\D_1$ is obtained from the second of Eqs.~\ref{stateequationevo} by substituting $t$ with the symbol $\D_1$, $\xd \oftime$ with the value of $\xd_s$ from Eq.~\ref{intersectionpoint_DI} and $u^*$ with the value of $u^*_1$  from Eq.~\ref{controloutofswitchcurve}. Finally, the expression of $\D_2$ is obtained from the same equation by substituting $t$ with the symbol $\D_2$, $\xd \oftime$ with zero, $u^*$ with the value of $u^*_2$  from Eq.~\ref{controloutofswitchcurve} and $\xd_0$ with the value of $\xd_s$ from Eq.~\ref{intersectionpoint_DI}.

Furthermore, the properties in Eqs.~\ref{DINR2Rlambineq},~\ref{DINR2Rlambineq33}, and \ref{DINR2RdiseqTstar} are proven by exhaustive analysis of validity for all possible sign combinations. See Tables~\ref{tableanalysisDINR2R1},~\ref{tableanalysislambdaeq2}  and~\ref{tableanalysisDINR2R2}, respectively.

Finally, by considering the following initial values for the costate
\be~\label{initialcostatecasef00}
\displaystyle p^*_{{\dot{x}}0} = p^*_{x0}\,\Delta_1, \qquad  p^*_{x0}=\frac{\rho^*}{\Lambda_0\Sigma_0}\, \qquad \forall \rho^*,
\ee
the zero-value condition for the Hamiltonian in Eq.~\ref{zerohamiltonianspecific} holds true, as it can be immediately verified. 
The values in Eq.~\ref{initialcostatecasef00} are obtained by considering the algebraic system formed by the equation obtained by imposing to zero the argument of the sign function in Eq.~\ref{regulatoroptimalcontrolexpl} with $t=t_s=\Delta_1$, and by the equation obtained by imposing to zero the Hamiltonian at the initial time (Eq.~\ref{hamiltonianforallt} with $t=0$).

Finally, Eqs~\ref{eq30} to ~\ref{eq32} are immediately demonstrated by integrating in time the first of Eq.~\ref{doubleintegrator2} between the points of interest and solving for the time. 
\begin{flushleft}
\hfill \bf Q.e.d.
\end{flushleft}
\renewcommand{\arraystretch}{1.4}

\begin{table}[t!]
\setlength{\tabcolsep}{9pt}
\begin{tabular}{@{}lllc@{}}
\toprule
 $\boldsymbol{\Sigma_0}$  & $\boldsymbol{\xd_0}$   & \thead{\bf{Eq.~\ref{DINR2Rlambineq} }\\\bf{yields}}&  \thead{\bf{Eq.~\ref{DINR2Rlambineq} }\\\bf{is}}\\
 \cmidrule(rl){3-3}
 \midrule
 \addlinespace
 1      &    $>0$  &$x_0>-\displaystyle \frac{\xd_0}{2} $  &  True*\\  
 \addlinespace
 1      &    $<0$  &$x_0>-\displaystyle \frac{\xd_0}{2} $  &  True*\\ 
 \addlinespace
-1      &    $>0$  &$x_0<\displaystyle \frac{\xd_0}{2} $  &  True*\\ 
\addlinespace
-1      &    $<0$  &$x_0<\displaystyle \frac{\xd_0}{2}  $ &  True*\\ 
\bottomrule
\end{tabular}
\vspace{0.05in}
\caption{\emph{Analysis of validity of Eq.~\ref{DINR2Rlambineq} }}
\label{tableanalysisDINR2R1}
\small{*[Because it agrees with Eq.~\ref{DINR2RF0}]}
\end{table}

\begin{table}[t!]
\setlength{\tabcolsep}{6pt}
\begin{tabular}{@{}llllc@{}}
\toprule
 $\boldsymbol{\Sigma_0}$  & $\boldsymbol{\xd_0}$   & \multicolumn{2}{c}{\bf{Eq.~\ref{definitionofF} yields}} & \thead{\bf{Eq.~\ref{DINR2Rlambineq33}}\\\bf{is}}\\
 \cmidrule(rl){3-4}
\midrule
\addlinespace
1   &   $<0$  & $\x_0 > \displaystyle \frac{\xd_0}{2}\rightarrow$ &$ \signFzero x_0 >\displaystyle \frac{\xd_0}{2}$ &True* \\ 
\addlinespace
  -1      &    $>0$  & $\x_0 <- \displaystyle \frac{\xd_0}{2}\rightarrow$ &$ \signFzero x_0 > \displaystyle \frac{\xd_0}{2}$ &True*  \\ 
\bottomrule
\end{tabular}
\vspace{0.05in}
\caption{\emph{Analysis of validity of Eq.~\ref{DINR2Rlambineq33}  }}
\label{tableanalysislambdaeq2}
\small{*[As obtained by considering together Eqs.~\ref{definitionofF} and~\ref{DINR2Rlambdef}]}
\end{table}

\begin{table}[t!]
\begin{tabular}{@{} l l l l l @{}}
\toprule
$\boldsymbol{\Sigma_0}$  & $\boldsymbol{\xd_0}$   &  \multicolumn{2}{c}{\bf{Eqs.~\ref{DINR2RTstarexplicit} and~\ref{DINR2RdiseqTstar} yields}}& \thead{\bf{Eq.~\ref{DINR2RdiseqTstar}}\\\bf{is}}\\
 \cmidrule(rl){3-4}
 \midrule
  1      &    $>0$  &$2\Lambda_0+\xd_0>\xd_0 \rightarrow$&$\Lambda_0 >0$  &  True*\\ 
\addlinespace
1      &    $<0$  &$2\Lambda_0-|\xd_0|>|\xd_0| \rightarrow$&$\Lambda_0 >|\xd_0|$  &  True**\\ 
\addlinespace
-1      &    $>0$  &$2\Lambda_0-|\xd_0|>|\xd_0| \rightarrow$&$\Lambda_0 >|\xd_0|$  &  True**\\ 
\addlinespace
-1      &    $<0$  &$2\Lambda_0-\xd_0>-\xd_0\rightarrow$&$ \Lambda_0 >0$  &  True*\\ 
\bottomrule
\end{tabular}
\vspace{0.05in}
\caption{\emph{Analysis of validity of Eq.~\ref{DINR2RdiseqTstar}}}
\label{tableanalysisDINR2R2}
\small{*[Because it agrees with Eq.~\ref{DINR2Rlambdeflt0}]\\
{**}[Because it agrees with Eq.~\ref{DINR2Rlambineq33}]}
\end{table}

\bibliographystyle{plain}        
\bibliography{bibliographydatabaseCurtiRomano}           

\end{document}